\documentclass[a4paper,11pt, reqno]{amsart}

\usepackage[margin=1in]{geometry}
\usepackage{hyperref}
\usepackage{amsmath,amsfonts}
\usepackage{amsthm}
\usepackage{amssymb}
\usepackage{url}
\def\R{\mathbb{R}}

\def\SPP{\mathrm{SPP}}
\usepackage{multirow}

\usepackage{adjustbox}

\newtheorem{theorem}{Theorem}

\newtheorem{example}[theorem]{Example}
\usepackage{enumitem}

\usepackage{setspace}
\usepackage[ruled,vlined]{algorithm2e}
\usepackage{tikz}

\usepackage{pgfplots}
\pgfplotsset{compat=newest}
\usetikzlibrary{decorations.markings}

\allowdisplaybreaks

\let\origmaketitle\maketitle
\def\maketitle{
  \begingroup
  \def\uppercasenonmath##1{} 
  \let\MakeUppercase\relax 
  \origmaketitle
  \endgroup
}

\title{}

\author{}

\title[]{\Large Network Flow based approaches for the\\Pipelines Routing Problem in Naval Design}

\author[V. Blanco, G. Gonz\'alez, Y. Hinojosa, D. Ponce, M.A. Pozo \MakeLowercase{and} J. Puerto]{{\large V\'ictor Blanco$^\dagger$, Gabriel Gonz\'alez$^\ddagger$, Yolanda Hinojosa$^\star$,\\ Diego Ponce$^\ddagger$,Miguel A. Pozo$^\ddagger$, Justo Puerto$^\ddagger$}\medskip\\
$^\dagger$Institute of Mathematics  (IMAG) and Dpt Quant. Methods Economics \& Business, Universidad de Granada\\
$^\star$ Institute of Mathematics (IMUS) and Dpt. Applied Economics I, Universidad de Sevilla\\
$^\ddagger$Institute of Mathematics (IMUS) and Dpt. Stats \& OR, Universidad de Sevilla
}

\date{\today}

\begin{document}

\maketitle

\begin{abstract}
In this paper we propose a general methodology for the optimal automatic routing of spatial pipelines motivated by a recent collaboration with  Ghenova, a leading Naval Engineering company. We provide a minimum cost multicommodity network flow based model for the problem incorporating all the technical requirements for a feasible pipeline routing. A branch-and-cut approach is designed and different matheuristic algorithms are derived for solving efficiently the problem. We report the results of a battery of computational experiments to assess the  problem  performance as well as a case study of a real-world naval instance provided by our partner company.
\keywords{Pipeline Routing, Network Design,  Branch-and-Cut, Matheuristics, Naval Engineering.}

\end{abstract}

\section{Introduction}

The design and construction of ships is a complex procedure that starts with a concept design phase. In the concept design, all the necessary ship construction requirements are determined, in order to meet the client needs. Next, in the basic design phase of the ship structure, piping and electrical circuits are defined and all the different  equipment are adequately selected. In this phase, one of the most difficult tasks is to determine the pipeline route design in which one  has to determine the paths of the different services attending to a series of technical requirements.  Pipeline routing requires the coordination of different elements of the ship, as the routing of electrical circuits, water pipes and optical fiber, the prevention of traversing forbidden obstacles, the assuredness of the adequate space for handling the machinery, among many others, all of them having to guarantee the compatibility and the manufacturability technical requirements. Although there are some available tools to help the designers on this task, the difficulty of the problem makes that an expert is still necessary to guide the design to a desired final product. Usually, there are some specifications regarding the number of main lines, branches to be connected to the main lines, valves, etcetera that have to be considered (see e.g.  \cite{Cuervas}). Once the initial technical requirements are determined, the pipeline routing problem must be solved following these specifications. The main limitations of the pipeline routing problem on a ship can be classified in three types (see \cite{and11,P02,QRW08}):
\begin{description}
\item[Physical Constraints:] The path followed by the pipelines must avoid physical obstacles and connect with the adequate equipment.
\item[Operational Constraints:] The routes must consider accessibility for handling equipment and valves slackness for security. In addition, some zones in the decision space are more desirable than others, as bottom spaces in a cabin that may be used for storing other types of materials.
\item[Economic Constraints:] There is a limited budget both for the material and workforce cost. Thus, it is necessary to reduce the pipe lengths, as well as the use of elbows in the design.
\end{description}
Some of the above mentioned requirements must be imposed to design feasible  routes, while others can be quantified and reflected in a cost function in order to evaluate the different feasible alternatives and decide the most favorable. Nevertheless, defining  {adequately} such a function implies including elements of different nature as  {material and workforce costs}, use of elbows, routes passing through preference zones and holes, etc.

There exist a few software tools (as AVEVA, FORAN or SMART3D, among others) that may help industrial designers for the development of this work.  However, even those tools that incorporate some automatic routing functionality are still limited to consider all the technical requirements of this real-world problem being them more adequate for industrial plants. Furthermore, the underlying algorithms of these modules do not allow to incorporate the different specificities for the design of a ship, and then, in most cases the solution  is not valid for the naval designer. In other cases, the decision aid tools are extremely specialized, not being flexible enough for the design of networks of different characteristics.

The goal of this paper is to derive an efficient mathematical programming based methodology for the automatic design of pipeline routes that take into account the different technical  {specificities} of a naval design, but still flexible enough to be adapted to different situations and constructibility criteria  {that appear in real-world situations as the one that motivated this study.}

Automatic pipeline routing has been already addressed in the literature (see for instance, \cite{gui05,Kim11,Min2020,par02,Singh2021}). One of the most crucial steps in this problem is the description of the set of feasible routes for the pipes. Although the problem is initially stated in a continuous framework, with routes that are allowed to be traced in a continuous three-dimensional region, the mathematical problem is intractable in this form (one would need to locate a set of complex three-dimensional paths on a continuous space). Thus, a discretization of the whole space into {a finite set of feasible routes} is needed also because the tractability of the problem and the simplicity of the obtained routes.  One of the possible options is to discretize the space by cell decomposition. In \cite{lee61} the author proposed to subdivide the region {into} squared cells and those containing obstacles are removed. Several authors have used that scheme to derive algorithms for solving shortest path problems in continuous regions~\cite{H69,mit86,rou75,sch99}. This approach has been also adapted to solve the pipeline routing problem \cite{and11,asm13,asm06,kim13,zhu91}. However, this method requires a high number of cells in the subdivision to get accurate solutions of the problem. To overcome this situation, one may also discretize the space by means of a graph  {structure}~\cite{gui05,new73,shi86}.  Both approaches, cell subdivision and graph generation can be also adequately combined to derive efficient approaches in pipeline routing~\cite{par02}.

Once the space is adequately discretized, there are several heuristic approaches that have been proposed to solve the pipeline routing problem. Most of them are genetic algorithms \cite{IKK05,Ito99,Kim11} or ant colony algorithms~\cite{XYZ06,XYZ07}.

In this paper we will consider a graph-based discretization of the space that allows us to use Combinatorial Optimization tools for {efficiently solving} the problem.  The graph is generated taking into account the shape of the region and obstacles, the positions of the source and destination points of the pipes, the size of the pipelines, the preferred zones and penetrable zones (holes through which some obstacles can be traversed), and other characteristics of the pipes to route. Firstly, we generate the nodes of the graph using highlighted points (source and destination points, corners, holes, ...) and also intermediate points by means of a grid with a desired width. Secondly, we define the edges of the graph by linking the close-enough nodes. Finally, each edge is provided with a set of weights (one for each pipeline to route) representing the different costs of using it in a path.

 On the basis of this idea our main contributions are:
\begin{enumerate}
\item We discretize the continuous three-dimensional space  {by} using a graph-based framework that allows searching for pipeline routes {taking} into account the obstacles and the use of elbows in the routes.
\item We propose a multicommodity network flow based model for the problem that incorporates the different physical and operational limitations of the routing. In particular, the minimum allowed distances between consecutive elbows and separation between services. We provide a particular branch-and-cut approach for the problem.
\item We propose a flexible assessment function to evaluate feasible routes that take into account the different specificities that appear in real-world naval design: length of the paths, preference zones, closeness to ceilings/floors, use of elbows, crossing penetrable zones,  etc.
\item We develop different matheuristic algorithms able to solve realistic instances of the problem in reasonable CPU times.
\end{enumerate}

To present our contribution we have organized the paper in six sections. Section \ref{sec:PRPND} describes the main elements involved in the pipeline routing problem and their mathematical representation.
Section \ref{sec:model} is devoted to present the mathematical programming model that we propose for solving the pipeline routing problem imposing naval design technical requirements.
In Section \ref{sec:exact} a branch-and-cut method is proposed in order to incorporate complicating constraints as they are required, instead of considering, initially, all of them. Two families of matheuristic algorithms are provided in Section \ref{sec:matheuristic} to solve large-sized instances. In Section \ref{sec:experiments} we report the results of  our computational experiments. There, we also include a case study based on real data provided by our industrial partner. Finally, we derive some conclusions and future lines of research is Section \ref{sec:conclusions}. We have included an Appendix to gather all the pseudocodes that describe the details of our algorithms.

\section{The Pipeline Routing Problem in Naval Design}\label{sec:PRPND}

The goal of this paper is to model and solve the problem of how to route different pipelines on a ship competing for a reduced space taking into account the adequate technical requirements.  In what follows we describe the main elements of the problem: input data, feasible actions and assessment of a particular solution.

We assume that the underlined region where the pipelines are to be routed is represented by a bounded polyhedron, $\mathcal{P}_0$, in $\R^3$. In real-world situations, the polyhedron is usually a cuboid in the form $\mathcal{P}_0 = [a_1,a_2] \times [b_1,b_2] \times [c_1,c_2]$, representing a cabin in the ship. We are also given a finite set of regions $O_1, \ldots, O_o \subseteq \mathcal{P}_0$ that represents $o$ obstacles that cannot be traversed by a pipeline route.  For the sake of simplicity, these  regions will be also identified with cuboids. Thus, the  three-dimensional region where the pipelines are allowed to be traced is $\mathcal{P} = \mathcal{P}_0 \backslash \cup_{\ell=1}^o O_\ell$.
In addition, we are given a finite set of services (cylindrical pipelines) each of them represented by the coordinates of a source and a  destination point (both in $\mathcal{P}$),   a radius (of the cylindrical pipeline), and a safety minimum allowed distance with respect to other services.  Our approach can be easily adapted to other shapes of pipelines, as parallelepipeds (with rectangular sections).

The main decision to be taken in this problem is to determine the set of (dimensional) paths in $\mathcal{P}$ that each one of the services must  follow.  Initially,  a feasible path for a service will be a continuous union of segments in $\mathcal{P}$ starting off from its source point and reaching its destination point. In addition, a minimum required distance must be ensured between two consecutive breakpoints (elbows) by {constructibility reasons.} Furthermore, once a feasible path for each service is traced,   and  once a dimensional  pipeline is traced around the path,  one has to ensure that a minimum required distance between   the different services is respected. In particular, the pipelines are not allowed to cross or to overlap in a solution.

The quality of a suitable solution of the problem, given by a set of feasible paths, is evaluated with respect to different criteria. Apart from the length, which is proportional to costs, one also looks for paths with small number of elbows, close to the ceiling of the region or traversing pre-specified preference zones.

In what follows we describe the mathematical representation of the elements involved in the Pipeline Routing Problem in a Ship (PRPS) that we analyze.

  In order to discretize   the continuous space $\mathcal{P}$, first, an orthogonal 3D grid is built on the big parallelepiped $\mathcal{P}_0$, assuring that source and destination points of  each service are nodes of the grid and that they are connected with other nodes of the grid. This grid defines a baseline undirected graph to which nodes and edges intersecting   regions $O_\ell$, for $\ell=1, \ldots, o$,  are removed. In Example \ref{ex:ilustra}  we  show a simple instace  for the problem.

  \begin{example} \label{ex:ilustra}
In order to illustrate the procedure described in the paper we include an example of a scenario with an obstacle and three services (see Figure \ref{fig:x} (left)). The scenario is a box with dimensions $6\times 6 \times 4$ with an interior obstacle (box) of dimensions $2\times 2 \times 4$ located in the center of the box representing a  pillar.  The  discretized space of solutions is given by a three-dimensional orthogonal grid from which  the corresponding portion that is  within the obstacle has been removed (see Figure \ref{fig:x} (right)).

\begin{figure}
\centering
\includegraphics[scale=0.45]{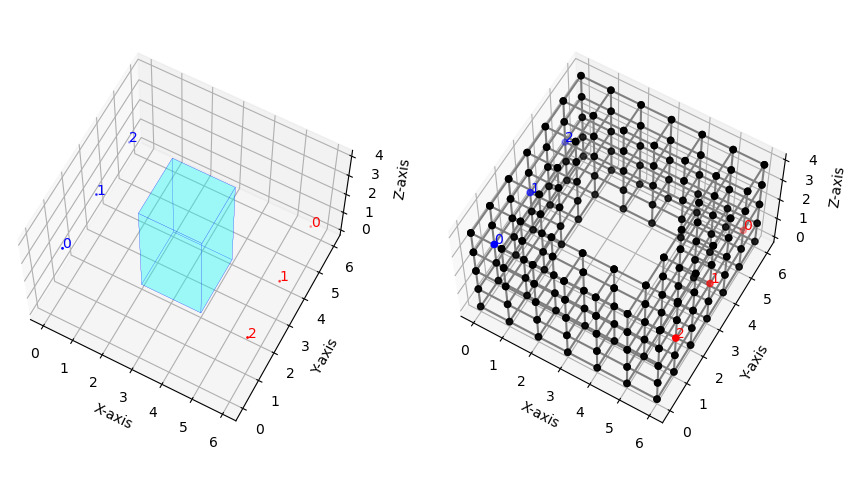}
\caption{Illustration of the first discretization phase of a PRPS  instance  with an obstacle and three services with sources and destinations at the blue and red points, respectively. The left figure shows the origin and destination points and the obstacle. The right figure depicts the grid of the discretized space.\label{fig:x}}
\end{figure}
\end{example}

   However,  when tracing a path in such a graph, one is not able to detect or penalize the use of elbows.  Elbows in a pipeline must be adequately identified in its route both because constructibility and also because routes with a smaller number of elbows are preferred. The goal is to obtain a  { mathematical programming model for the problem that avoids non-linearities.  In order to consider elbows in a linear objective function in our model,} we modify the initial graph by \textit{exploiting} the nodes of the graph as follows:
\begin{enumerate}
\item Each physical node $v$ in the initial graph is replaced by an \textit{exploited node}, i.e., a set of three \textit{virtual nodes}, $v_X, v_Y, v_Z$ with the same 3D coordinates than $v$, one for each direction of the canonical basis of $\R^3 $ (X, Y and Z).
\item Virtual nodes associated to the same physical node are linked through the so-called \textit{virtual edges}. Their lengths are zero (since they link nodes with the same coordinates), but they will have a positive cost representing the usage of an elbow in the route.
\item Each physical edge in the initial graph is replaced by another edge with the same length as the physical one, where instead of joining two physical nodes, it connects two virtual nodes (associated to the same physical end nodes of the original edge). Specifically, each edge is linked to the virtual nodes identified with its direction. In this way, edges parallel to the X-axis link $X$-virtual nodes, edges parallel to the Y-axis link $Y$-virtual nodes, and edges parallel to the Z-axis link $Z$-virtual nodes.

\end{enumerate}

In Figure \ref{fig4} we show an illustration of this explosion of nodes. In the left picture we show a crossing node, $v$, which is linked in the graph to other six nodes. In the right picture we show the three exploited nodes, $v_X$, $v_Y$ and $v_Z$, of node $v$. Node $v_X$ (resp. $v_Y$, $v_Z$) is linked only to the two adjacent nodes which are linked with edges parallel to the $X$-axis (resp. $Y$-axis, $Z$-axis) and with the others two virtual nodes associated with the same physical node. As we can see,  the three virtual edges linking $v_X$ with $v_Y$, $v_X$ with $v_Z$ and $v_Y$ with $v_Z$ are drawn with dashed lines in this picture.

\usetikzlibrary{calc,3d}

\begin{figure}[h]
\begin{center}
\begin{tikzpicture}[scale=2]

	\coordinate (P1) at (-6cm,1.5cm); 
	\coordinate (P2) at (5cm,1.5cm); 

	\coordinate (A1) at (0em,0cm); 
	\coordinate (A2) at (0em,-2cm); 

	\coordinate (A3) at ($(P1)!.8!(A2)$); 
	\coordinate (A4) at ($(P1)!.8!(A1)$);

	\coordinate (A7) at ($(P2)!.7!(A2)$);
	\coordinate (A8) at ($(P2)!.7!(A1)$);

	\coordinate (A5) at
	  (intersection cs: first line={(A8) -- (P1)},
			    second line={(A4) -- (P2)});
	\coordinate (A6) at
	  (intersection cs: first line={(A7) -- (P1)},
			    second line={(A3) -- (P2)});

	\draw[thick,black] (A5) -- (A6);
	\draw[thick,gray] (A3) -- (A6);
	\draw[thick,gray!40] (A7) -- (A6);

	
	\coordinate (A9) at ($(A6) + (A6)-(A3)$);
	\draw[thick,gray] (A6) -- (A9);
	
	\coordinate (A10) at ($(A6) + (A6)-(A5)$);
	\draw[thick,black] (A6) -- (A10);

	\coordinate (A11) at ($(A6) + (A6)-(A7)$);
	\draw[thick,gray!40] (A6) -- (A11);

	\draw[fill=gray] (A3) circle (0.15em);
	\draw[fill=gray] (A9) circle (0.15em);
	\draw[fill=gray] (A5) circle (0.15em);
	\draw[fill=gray] (A10) circle (0.15em);
	\draw[fill=gray] (A7) circle (0.15em);
	\draw[fill=gray] (A11) circle (0.15em);
	
	\draw[fill=gray!30, inner sep=3mm] (A6) circle (0.4em) node {$v$};

\end{tikzpicture}~\begin{tikzpicture}[scale=2]

	\coordinate (P1) at (-6cm,1.5cm); 
	\coordinate (P2) at (5cm,1.5cm); 

	\coordinate (A1) at (0em,0cm); 
	\coordinate (A2) at (0em,-2cm); 

	\coordinate (A3) at ($(P1)!.8!(A2)$); 
	\coordinate (A4) at ($(P1)!.8!(A1)$);

	\coordinate (A7) at ($(P2)!.7!(A2)$);
	\coordinate (A8) at ($(P2)!.7!(A1)$);

	\coordinate (A5) at
	  (intersection cs: first line={(A8) -- (P1)},
			    second line={(A4) -- (P2)});
	\coordinate (A6) at
	  (intersection cs: first line={(A7) -- (P1)},
			    second line={(A3) -- (P2)});

	\draw[thick,black] (A5) -- (A6);
	\draw[thick,gray] (A3) -- (A6);
	\draw[thick,gray!40] (A7) -- (A6);

	
	\coordinate (A9) at ($(A6) + (A6)-(A3)$);
	\draw[thick,gray] (A6) -- (A9);
	
	\coordinate (A10) at ($(A6) + (A6)-(A5)$);
	\draw[thick,black] (A6) -- (A10);

	\coordinate (A11) at ($(A6) + (A6)-(A7)$);
	\draw[thick,gray!40] (A6) -- (A11);

	\draw[fill=gray] (A3) circle (0.15em);
	\draw[fill=gray] (A9) circle (0.15em);
	\draw[fill=gray] (A5) circle (0.15em);
	\draw[fill=gray] (A10) circle (0.15em);
	\draw[fill=gray] (A7) circle (0.15em);
	\draw[fill=gray] (A11) circle (0.15em);
	
	
	\coordinate (VX) at ($(A6) + 0.2*(A7)-0.2*(A6)$);
	\coordinate (VY) at ($(A6) + 0.2*(A5)-0.2*(A6)$);
	 \coordinate (VZ) at ($(A6) + 0.2*(A3)-0.2*(A6)$);

	\draw[fill=gray!30, inner sep=1mm] (VX) circle (0.3em) node {\tiny $v_X$};
	\draw[fill=gray!30, inner sep=1mm] (VY) circle (0.3em) node {\tiny $v_Z$};
	\draw[fill=gray!30, inner sep=1mm] (VZ) circle (0.3em) node {\tiny $v_Y$};

	\draw[dashed] (VX)--(VY);
	\draw[dashed] (VX)--(VZ);
	\draw[dashed] (VZ)--(VY);

	\draw[fill=gray!30, inner sep=1mm] (VX) circle (0.3em) node {\tiny $v_X$};
	\draw[fill=gray!30, inner sep=1mm] (VY) circle (0.3em) node {\tiny $v_Z$};
	\draw[fill=gray!30, inner sep=1mm] (VZ) circle (0.3em) node {\tiny $v_Y$};

\end{tikzpicture}
\end{center}
\caption{Illustration of exploited nodes and virtual edges\label{fig4}}
\end{figure}
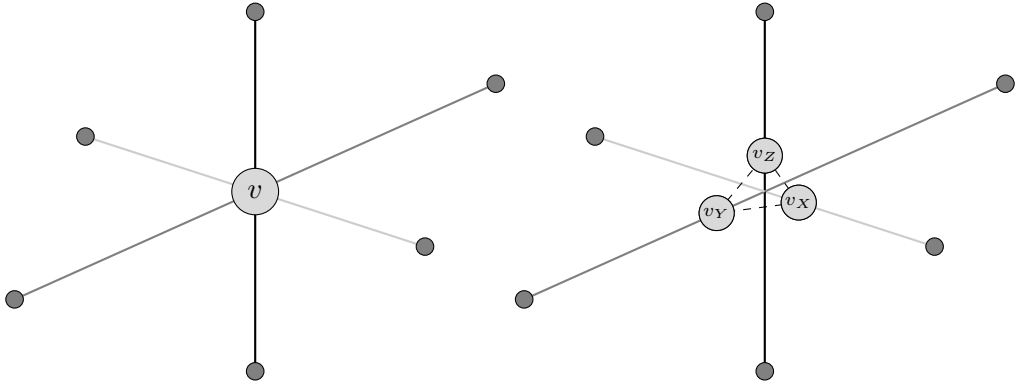

\begin{figure}
\begin{center}
\begin{tikzpicture}[scale=1.5]

	\coordinate (P1) at (-6cm,1.5cm); 
	\coordinate (P2) at (5cm,1.5cm); 

	\coordinate (A1) at (0em,0cm); 
	\coordinate (A2) at (0em,-2cm); 

	\coordinate (A3) at ($(P1)!.8!(A2)$); 
	\coordinate (A4) at ($(P1)!.8!(A1)$);

	\coordinate (A7) at ($(P2)!.7!(A2)$);
	\coordinate (A8) at ($(P2)!.7!(A1)$);

	\coordinate (A5) at
	  (intersection cs: first line={(A8) -- (P1)},
			    second line={(A4) -- (P2)});
	\coordinate (A6) at
	  (intersection cs: first line={(A7) -- (P1)},
			    second line={(A3) -- (P2)});

	\draw[thin,black] (A5) -- (A6);
	\draw[thin,gray] (A3) -- (A6);
	\draw[thin,gray!40] (A7) -- (A6);

	
	\coordinate (A9) at ($(A6) + (A6)-(A3)$);
	\draw[thin,gray] (A6) -- (A9);
	
	\coordinate (A10) at ($(A6) + (A6)-(A5)$);
	\draw[thin,black] (A6) -- (A10);

	\coordinate (A11) at ($(A6) + (A6)-(A7)$);
	\draw[thin,gray!40] (A6) -- (A11);

	\draw[fill=gray!5] (A3) circle (0.15em) node {\tiny $a$};
	\draw[fill=gray!5] (A9) circle (0.15em) node {\tiny $b$};
	\draw[fill=gray] (A5) circle (0.15em);
	\draw[fill=gray] (A10) circle (0.15em);
	\draw[fill=gray] (A7) circle (0.15em);
	\draw[fill=gray] (A11) circle (0.15em);
	
	
	\coordinate (VX) at ($(A6) + 0.2*(A7)-0.2*(A6)$);
	\coordinate (VY) at ($(A6) + 0.2*(A5)-0.2*(A6)$);
	 \coordinate (VZ) at ($(A6) + 0.2*(A3)-0.2*(A6)$);
	
	\node[circle,draw,fill=gray!30,inner sep=0.05mm,minimum size=1mm] (VVZ) at  (VZ) {\tiny $v_Y$};
	\node[circle,draw,fill=gray!30,inner sep=0.05mm,minimum size=1mm] (VVY) at  (VY) {\tiny $v_Z$};
		\node[circle,draw,fill=gray!30,inner sep=0.05mm,minimum size=1mm] (VVX) at  (VX) {\tiny $v_X$};

	\draw[dashed] (VVX)--(VVY);
	\draw[dashed] (VVX)--(VVZ);
	\draw[dashed] (VVZ)--(VVY);
	
	\draw[very thick] (A3)--(VVZ);
	\draw[very thick] (VVZ)--(A9);

	\node[circle,draw,fill=gray!5,inner sep=0.05mm,minimum size=3mm] () at  (A3) {\scriptsize$a$};
	\node[circle,draw,fill=gray!5,inner sep=0.05mm,minimum size=3mm] () at  (A9) {\scriptsize$b$};
\end{tikzpicture}~\begin{tikzpicture}[scale=1.5]

	\coordinate (P1) at (-6cm,1.5cm); 
	\coordinate (P2) at (5cm,1.5cm); 

	\coordinate (A1) at (0em,0cm); 
	\coordinate (A2) at (0em,-2cm); 

	\coordinate (A3) at ($(P1)!.8!(A2)$); 
	\coordinate (A4) at ($(P1)!.8!(A1)$);

	\coordinate (A7) at ($(P2)!.7!(A2)$);
	\coordinate (A8) at ($(P2)!.7!(A1)$);

	\coordinate (A5) at
	  (intersection cs: first line={(A8) -- (P1)},
			    second line={(A4) -- (P2)});
	\coordinate (A6) at
	  (intersection cs: first line={(A7) -- (P1)},
			    second line={(A3) -- (P2)});

	\draw[thin,black] (A5) -- (A6);
	\draw[thin,gray] (A3) -- (A6);
	\draw[thin,gray!40] (A7) -- (A6);

	
	\coordinate (A9) at ($(A6) + (A6)-(A3)$);
	\draw[thin,gray] (A6) -- (A9);
	
	\coordinate (A10) at ($(A6) + (A6)-(A5)$);
	\draw[thin,black] (A6) -- (A10);

	\coordinate (A11) at ($(A6) + (A6)-(A7)$);
	\draw[thin,gray!40] (A6) -- (A11);

	\draw[fill=gray!5] (A3) circle (0.15em) node {\tiny $a$};
	\draw[fill=gray] (A9) circle (0.15em);
	\draw[fill=gray!5] (A5) circle (0.15em)  node {\tiny $b$};
	\draw[fill=gray] (A10) circle (0.15em);
	\draw[fill=gray] (A7) circle (0.15em);
	\draw[fill=gray] (A11) circle (0.15em);
	
	
	\coordinate (VX) at ($(A6) + 0.2*(A7)-0.2*(A6)$);
	\coordinate (VY) at ($(A6) + 0.2*(A5)-0.2*(A6)$);
	 \coordinate (VZ) at ($(A6) + 0.2*(A3)-0.2*(A6)$);

	\node[circle,draw,fill=gray!30,inner sep=0.05mm,minimum size=1mm] (VVZ) at  (VZ) {\tiny $v_Y$};
	\node[circle,draw,fill=gray!30,inner sep=0.05mm,minimum size=1mm] (VVY) at  (VY) {\tiny $v_Z$};
		\node[circle,draw,fill=gray!30,inner sep=0.05mm,minimum size=1mm] (VVX) at  (VX) {\tiny $v_X$};

	\draw[dashed] (VVX)--(VVY);
	\draw[dashed] (VVX)--(VVZ);
	\draw[dashed] (VVZ)--(VVY);
	
	\draw[very thick] (A3)--(VVZ);
	\draw[very thick] (VVZ)--(VVY);
	\draw[very thick] (VVY)--(A5);

	\node[circle,draw,fill=gray!5,inner sep=0.05mm,minimum size=3mm] () at  (A3) {\scriptsize$a$};
	\node[circle,draw,fill=gray!5,inner sep=0.05mm,minimum size=3mm] () at  (A5) {\scriptsize$b$};
\end{tikzpicture}~\begin{tikzpicture}[scale=1.5]

	\coordinate (P1) at (-6cm,1.5cm); 
	\coordinate (P2) at (5cm,1.5cm); 

	\coordinate (A1) at (0em,0cm); 
	\coordinate (A2) at (0em,-2cm); 

	\coordinate (A3) at ($(P1)!.8!(A2)$); 
	\coordinate (A4) at ($(P1)!.8!(A1)$);

	\coordinate (A7) at ($(P2)!.7!(A2)$);
	\coordinate (A8) at ($(P2)!.7!(A1)$);

	\coordinate (A5) at
	  (intersection cs: first line={(A8) -- (P1)},
			    second line={(A4) -- (P2)});
	\coordinate (A6) at
	  (intersection cs: first line={(A7) -- (P1)},
			    second line={(A3) -- (P2)});

	\draw[thin,black] (A5) -- (A6);
	\draw[thin,gray] (A3) -- (A6);
	\draw[thin,gray!40] (A7) -- (A6);

	
	\coordinate (A9) at ($(A6) + (A6)-(A3)$);
	\draw[thin,gray] (A6) -- (A9);
	
	\coordinate (A10) at ($(A6) + (A6)-(A5)$);
	\draw[thin,black] (A6) -- (A10);

	\coordinate (A11) at ($(A6) + (A6)-(A7)$);
	\draw[thin,gray!40] (A6) -- (A11);

	\draw[fill=gray!5] (A3) circle (0.15em) node {\tiny $a$};
	\draw[fill=gray] (A9) circle (0.15em);
	\draw[fill=gray] (A5) circle (0.15em);
	\draw[fill=gray] (A10) circle (0.15em);
	\draw[fill=gray] (A7) circle (0.15em);
	\draw[fill=gray!5] (A11) circle (0.15em) node {\tiny $b$};
	
	
	\coordinate (VX) at ($(A6) + 0.2*(A7)-0.2*(A6)$);
	\coordinate (VY) at ($(A6) + 0.2*(A5)-0.2*(A6)$);
	 \coordinate (VZ) at ($(A6) + 0.2*(A3)-0.2*(A6)$);

	\node[circle,draw,fill=gray!30,inner sep=0.05mm,minimum size=1mm] (VVZ) at  (VZ) {\tiny $v_Y$};
	\node[circle,draw,fill=gray!30,inner sep=0.05mm,minimum size=1mm] (VVY) at  (VY) {\tiny $v_Z$};
		\node[circle,draw,fill=gray!30,inner sep=0.05mm,minimum size=1mm] (VVX) at  (VX) {\tiny $v_X$};

	\draw[dashed] (VVX)--(VVY);
	\draw[dashed] (VVX)--(VVZ);
	\draw[dashed] (VVZ)--(VVY);
	
	\draw[very thick] (A3)--(VVZ);
	\draw[very thick] (VVZ)--(VVX);
	\draw[very thick] (VVX)--(A11);

	\node[circle,draw,fill=gray!5,inner sep=0.05mm,minimum size=3mm] () at  (A3) {\scriptsize$a$};
	\node[circle,draw,fill=gray!5,inner sep=0.05mm,minimum size=3mm] () at  (A11) {\scriptsize$b$};
\end{tikzpicture}
\end{center}
\caption{Illustration of the use of exploited nodes and virtual edges to model different types of turns in a path from node $a$ to node $b$\label{fig:elbows} }
\end{figure}
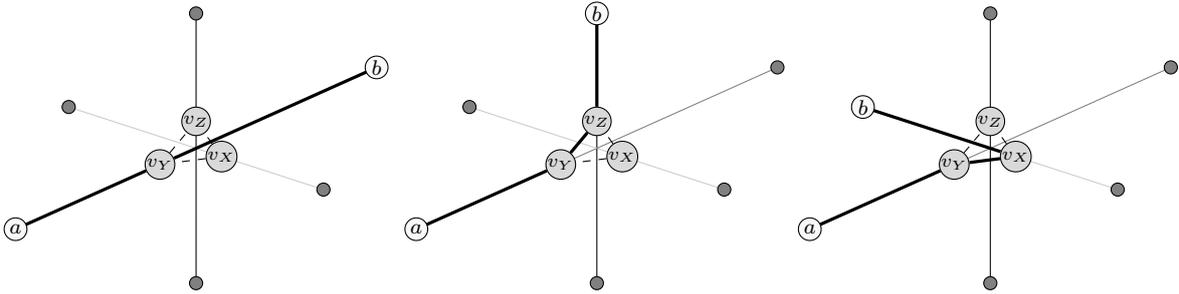

In Figure \ref{fig:elbows} we illustrate the use of this graph to model different types of turns in a path. In the left picture, to go from $a$ to $b$ no elbows are used and a single exploited node ($v_Y$) is used. On the other hand, in the center (resp. right) picture, turning from an edge parallel to $Y$-axis towards an edge parallel to $Z$-axis (resp. $X$-axis) requires traversing the virtual edge $(v_Y,v_Z)$ (resp. $(v_Y,v_X)$) incurring in a cost for using that elbow. The result is a graph representing the feasible space where routing the pipelines.

We denote by $\widetilde{G}=(V,E)$ the graph constructed as described above where $V$ is the set of (virtual) nodes  and $E$ the set of edges. Set $E$ includes the set of virtual edges, denoted by $E^v$, and  the replicas of the physical initial edges.  Since $\widetilde{G}$ is embedded in $\R^3$, each node $v \in V$ is also identified with its $3D$ coordinates, $(x_v, y_v, z_v) \in \R^3$. We denote by $d_{v v^\prime} = d((x_v,y_v,z_v), (x_{v^\prime},y_{v^\prime},z_{v^\prime}))$ the Euclidean distance between  nodes $v, {v^\prime} \in V$ and by $d_{ee^\prime}$ the Euclidean distance between the edges $e, e^\prime \in E$ (representing either, a point in $\R^3$ if $e \in E^v$ or a segment in $\R^3$ if $e \in E\setminus E^v$).

We assume that we are given a finite set of commodities, $K$, each of them requiring  its own pipeline to determine the design of its route. Each commodity $k\in K$ is defined by a pair $(s^k, t^k)$ where $s^k \in V$ is the source of the commodity and $t^k\in V$ is the destination of the commodity. Depending on the commodity $k\in K$, a cost structure, $c_e^k\geq 0$ for each $e \in E$,  is defined over the edges of the graph. These costs depend on many factors as for instance, the length of the edges, the elbow costs, the preferences  of the designer when routing the pipelines, etc. A detailed cost structure will be described in  {Section \ref{sec:case}} for the case study we deal with in this work.

\section{A multicommodity network flow based model}\label{sec:model}

In this section we describe the mathematical programming model that we propose for solving the PRPS imposing naval design technical requirements.  In a first approximation, one can model the problem as a minimum cost multicommodity network flow problem (\mbox{MCMNFP}, for short). Let us denote by $G=(V,A)$ (where $A=\{(i,j)\cup (j,i): \ e=\{i,j\} \in E\}$ is the arc set) the directed version of the graph $\widetilde{G}$ described in the previous section. Similar to the non-directed version, we assume that set $A$ includes the set of virtual arcs, denoted by $A^v$, and the replicas of the physical arcs. We also denote by $d_{aa^\prime}$ the Euclidean distance between the arcs $a, a^\prime \in A$ (notice that $d_{aa^\prime} =0$ if $a$ and $a^\prime$ are the arcs corresponding to the two opposite directions of the same edge). In addition, we assume that, for each commodity $k\in K$, the cost of each arc  $(i,j)\in A$ is
$ c_{ij}^k=c_{ji}^k= c_e^k$, where $e=\{i,j\}$.

The goal of MCMNFP is to route jointly  all the commodities in $K$ from their respective sources, $s^k$, to their destinations, $t^k$, through the network $G$ at minimum cost. This model has been widely studied in the literature ({the interested reader is referred to  \cite{Sal2020} and the references therein for further details on the MCMNFP}).

In order to describe the mathematical programming model for the MCMNFP, we use the following set of binary decision variables:
$$
x_{ij}^k = \left\{\begin{array}{cl}
1 & \mbox{if arc $(i,j)$ is used by the route of commodity $k$},\\
0 & \mbox{otherwise.}
\end{array}\right.
$$
With this set of variables, the overall cost of routing the commodities through the network can be expressed as:
$$
\sum_{(i,j) \in A}  \sum_{k\in K} c^k_{ij} x^k_{ij}.
$$
The following  set of linear constraints ensures the correct representation of the problem:
\begin{subequations}
    \makeatletter
        \def\@currentlabel{${\rm MCMNFP}$}
        \makeatother
       \label{MCMNFP}
        \renewcommand{\theequation}{${\rm MCMNFP}_{\arabic{equation}}$}
\begin{align}
&\sum_{k\in K} \sum_{j\in V: (i,j)\in A} x_{ij}^k  \leq 1, & \forall i \in V,\label{mcf:1}\\
&\sum_{j \in V: (i,j)\in A} x_{ij}^k - \sum_{j \in V: (i,j)\in A} x_{ji}^k = 0, & \forall k \in K, i \in V (i \neq s^k, t^k), \label{mcf:2}\\
&\sum_{j \in V: (s^k,j)\in A} x^k_{s^kj} = 1, &\forall k \in K,\label{mcf:3}\\
&\sum_{j \in V: (j,t^k)\in A} x^k_{jt^k} = 1, & \forall k \in K,\label{mcf:4}\\
&x_{ij}^k \in \{0,1\}, & \forall (i,j)\in A, k \in K.\label{mcf:5}
\end{align}
\end{subequations}
where \eqref{mcf:1} assures that only one commodity is routed through  a node and in particular, through an arc, \eqref{mcf:2}-\eqref{mcf:4} are the flow conservation constraints and \eqref{mcf:5} is the domain of the variables.

\begin{figure}[h]
\centering
\includegraphics[scale=0.6]{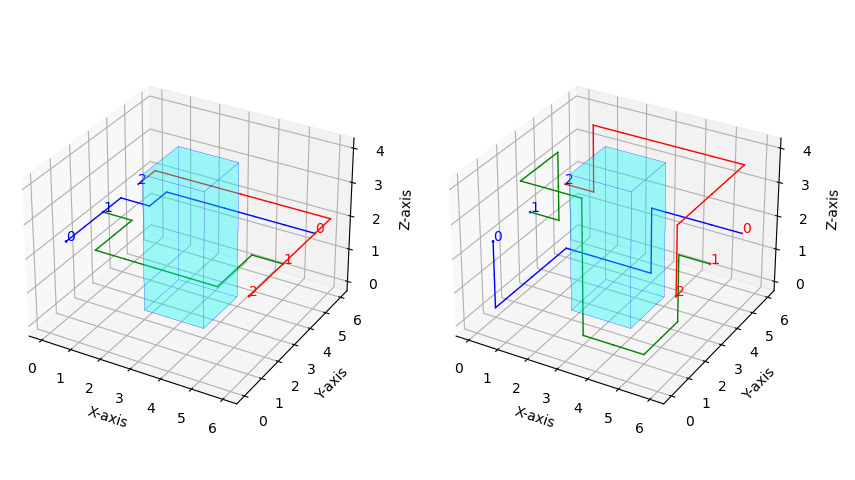}
\caption{Solution of a MCMNFP in the underlined graph without additional modelling constraints (left) and the final solution after adding the additional constraints on distances between services and elbow test  (right) \label{fig:xx}}
\end{figure}

The solution of MCMNFP is a set of  paths (one for each commodity) connecting sources with destinations  that do not overlap on the graph $G$. The reader may observe that in the graph $G$  each physical node has been replaced by three virtual nodes on the same coordinates. This explains that in Figure \ref{fig:xx} (left) some paths occupy the same space although they do not actually overlap on $G$.

This problem is known to be $\mathcal{NP}$-hard, even for the case of two commodities (see \cite{even1975complexity}). Additionally,  in pipeline routing there are two main technical requirement that must be considered in order to construct feasible designs: a minimum distance between the pipelines and a minimum distance between consecutive elbows in a single pipeline. The first requirement allows to adequately separate the different services to be used by the pipelines (and such that, when they are enlarged to represent dimensional pipes, they fit on the space)  as well as to keep space for the manipulation of machinery in the ship cabins. The second one, is a physical limitation of pipeline routing since there is not enough space for a pipe to turn  twice in case consecutive elbows are positioned too close. Note that solutions consists of segments in $\R^3$ with no dimensionality and therefore, these considerations are not taken into account in the MCMNFP. Thus, to impose those technical requirements the following constraints must be added to the problem in order to adequately model the PRPS:

\begin{enumerate}
\item {\it Distance between Services:}

Given  {an (non-virtual) active edge for a service $k$, a minimum distance with respect to} other services is required. Let $R^{k}$ be 
the radius of the cylinder of pipeline $k\in K$ and let $\Delta^k$ be the minimum security distance from $k$ to any other element. Then, the minimum distance allowed between services $k$ and $k^\prime$ is $R^{kk^\prime} = R^k+R^{k'}+\max\{\Delta^k,\Delta^{k^\prime}\}$.

This requirement is assured imposing the following set of constraints:
\begin{align}
	 & \sum_{k'\in K: k'\neq k}\sum_{\substack{a'\in A\\ d_{aa'}< R^{kk'}}}x_{a'}^{k'} \leq M_a^k(1-x_{a}^k),   & \forall a\in A\setminus A^v, k\in K, \label{mcf:6}\tag{${\rm Dist}$}
\end{align}
that is, if $a$ is an active arc in the path of service $k$, then, no arc in the path of any other service $k'\neq k$ can be activated at a distance smaller than the minimum required. Here, fixed $a \in A\setminus A^v$ and $k\in K$,  $M_a^k$ is a big enough constant (greater than $2(|K|-1)$ times the number of the grid edges contained in the cuboid centered at midpoint of arc $a$ and with edge length $length(a)+\, 2 \max_{k'\neq k} R^{kk'}$).

\item {\it Elbow Test}:

A feasible path for a service, $k \in K$, is required to verify a minimum distance between consecutive elbows, $D^k \geq 0$, of this single service. Since elbows are identified in the non-directed graph with virtual edges and then, with virtual arcs in the directed graph, this requirement can be incorporated to the model as follows:
\begin{align}
	 &  x_a^k+x_{a'}^k\leq 1,  &  \forall a, a'\in A^v: d_{aa'}\leq D^k, a\neq a', k\in K,\label{mcf:7}\tag{${\rm Elbow}$}
\end{align}
that is, if two different virtual arcs  (elbows) of the same service $k\in K$ are at a distance smaller than or equal to $D^k$, then only one of them is allowed to be active.

Example \ref{ex:ilustra2} shows the effects of considering constraints \eqref{mcf:6} and \eqref{mcf:7}.
\begin{example}[Example \ref{ex:ilustra} - continuation] \label{ex:ilustra2}

Let $(N,E)$  be the grid depicted  in Figure \ref{fig:x} (right). Grid $(N,E)$ is a simplification of the actual graph $G=(V,A)$ used to solve the problem because in the final graph  a directed version is considered and each of the nodes in the grid is replaced by its explosion as explained in Section \ref{sec:PRPND}. 

We have solved the Example \ref{ex:ilustra} scenario with the multicommodity model with and without the additional modelling constraints  \eqref{mcf:6} and \eqref{mcf:7}.  Observe that the paths without additional modelling constraints do not respect the conditions about distances and thus a much shorter solution is provided (see Figure \ref{fig:xx} (left)) than in the case in which the technical requirements imposed by means of constraints \eqref{mcf:6} and \eqref{mcf:7} have to be fulfilled  (see Figure \ref{fig:xx} (right)).

 Figure \ref{fig:finalsol} shows the final solution with dimensional pipelines of the given pre-specified radii.
\begin{figure}[h]
\centering
\includegraphics[scale=0.85]{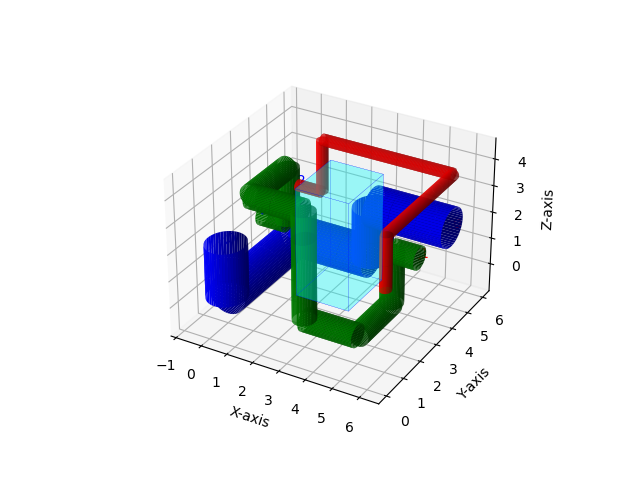}
\caption{ A graphical display of the solution of Example \ref{ex:ilustra} with the final pipelines of the three services and their actual radii \label{fig:finalsol}}
\end{figure}
\end{example}

Note that constraints \eqref{mcf:7} can be reinforced as follows. Given a virtual arc $a\in A^v$, we denote by $RE(a)=\{\tilde a\in A^v: \tilde a \mbox{ has the same coordinates as } a\}$, namely the set of arcs linking virtual nodes associated to the same physical node as the one linked by $a$. Clearly, if $d_{aa'}\leq D^k$, all the arcs in $RE(a)$ and $RE(a')$ are incompatible and then cannot take simultaneously a value of $1$. Thus,

\begin{align}
	 & \sum_{\tilde a \in RE(a)} x_{\tilde a}^k+ \sum_{\bar a \in RE(a')} x_{\bar a}^k\leq 1,  &  \forall a, a'\in A^v: d_{aa'}\leq D^k, a\neq a', k\in K.\label{mcf:8}\tag{${\rm Elbow-R}$}
\end{align}
\end{enumerate}

With the above considerations, the PRPS can be formulated as:
\begin{align}
\min & \sum_{(i,j) \in A}  \sum_{k\in K} c^k_{ij} x^k_{ij} \label{prps}\tag{${\rm PRPS}$}\\
\mbox{s.t. } & \eqref{mcf:1}-{\eqref{mcf:5}},\nonumber\\
& \eqref{mcf:6}, {\eqref{mcf:8}}.\nonumber
\end{align}

Problem \eqref{prps} is a minimum cost multicommodity network flow problem with additional constraints that enforce the fulfillment of the technical requirements for the pipelines in the solution. Although the problem is modeled as an Integer Linear Programming Problem, in real-world situations it has a large number of variables ($\mathcal{O}(|A||K|)$) and a large number of constraints. In particular, the technical constraints \eqref{mcf:6} and \eqref{mcf:8} place a considerable load on the model ($\mathcal{O}(|A|^2 |K|)$).

\subsection{An exact solution method for \eqref{prps}}\label{sec:exact}

The large dimensions of  our model makes the straightforward approach of putting it into a MIP solver not possible even for small-sized instances. We overcome this drawback using a branch-and-cut method that incorporates the constraints as they are required, instead of considering initially all of them.
This strategy allows an efficient exact solution approach of the problem  by means of solving an incomplete (relaxed) formulation with only some of the constraints in the model,  while the remaining constraints, required to assure the feasibility of the solutions, are incorporated \textit{on-the-fly} as needed. Although in the worst case situation, the procedure may need all the constraints  not initially included, in practice, only a small number of them are added, reducing considerably the size of the problem.

Specifically, we consider, in the beginning, the  relaxed master problem \eqref{MCMNFP}:
\begin{align}
\min & \sum_{(i,j) \in A}  \sum_{k\in K} c^k_{ij} x^k_{ij} \label{master}\tag{${\rm MCMNFP}$}\\
\mbox{s.t. } & \eqref{mcf:1}-{\eqref{mcf:5}}.\nonumber
\end{align}

The above problem is nothing but the multicommodity network flow model that does not consider the two families of technical constraints that are required for a feasible solution of our problem. Thus, when solving \eqref{master}, one may obtain solutions which do not satisfy constraints \eqref{mcf:6} and/or  \eqref{mcf:8}. Therefore, each obtained solution  must be checked for feasibility. To separate the violated constraints, we apply an enumerative procedure and those which are violated are added to the pool so that the problem is solved again (cutting off the previously obtained solution).

 {More specifically, for a given solution of the relaxed master problem, we check whether it violates \eqref{mcf:6} and/or  {\eqref{mcf:8} } as follows:
\begin{itemize}
\item To check the violation of any of the constraints in \eqref{mcf:6} we proceed by  measuring the distance between two arcs from different commodities in the solution or the distance between an arc from a commodity in the solution and obstacles' edges. If the distance measured is smaller than the sum of the radius of both commodities plus the security distance between them, the constraint is violated. In case it is violated for a service $k\in K$ and $a \in A \setminus A^v$, the constraint is added to the constraint pool of the master problem.

\item The violation of the constraints  {\eqref{mcf:8} } is tested by  first ordering virtual arcs within  each commodity in the path given by the current solution and then measuring the distance between two consecutive virtual arcs in the path. If the distance measured is smaller than the minimum required distance between elbows  for this commodity then the constraint is violated.  In this case, this particular constraint (in its reinforced version) is introduced  as a feasibility cut to the set of constraints of the problem.
\end{itemize}

For a more efficient implementation of the procedure, it is embedded in the branch-and-bound tree by means of lazy cuts, resulting in a branch-and-cut strategy for solving the problem. }

\section{ Two matheuristic algorithms for the  \ref{prps}}\label{sec:matheuristic}

 {In order to solve the pipeline routing model \ref{prps} presented in Section \ref{sec:model} one has to solve a MILP problem. However, in most cases of real-world situations, due to their large number of variables and constraints, the corresponding problem cannot be solved by a commercial solver  (and in many cases it cannot be even loaded). For this reason, we provide in this section  two families of matheuristic algorithms for the problem based on two different paradigms: reducing the dimension of the MILP model and decomposing the problem into simpler problems. The pseudocode of both algorithms is included in
\ref{appendix}.}

\subsection{Dimensionality Reduction\label{ss:dm}}

Recall that we assume that the underlined region where the pipelines are to be routed is represented by a bounded polyhedron in $\R^3$, and that  in real-world situations, the polyhedron is usually a cuboid representing a cabin in the ship where some other smaller cuboids (obstacles) have been removed. In addition, there are zones of such a region that are rarely used by the paths in the solutions (by all or some of the services). We exploit this fact  in order to reduce the number of variables of our model.

\begin{figure}[h]
\centering
\includegraphics[scale=0.7]{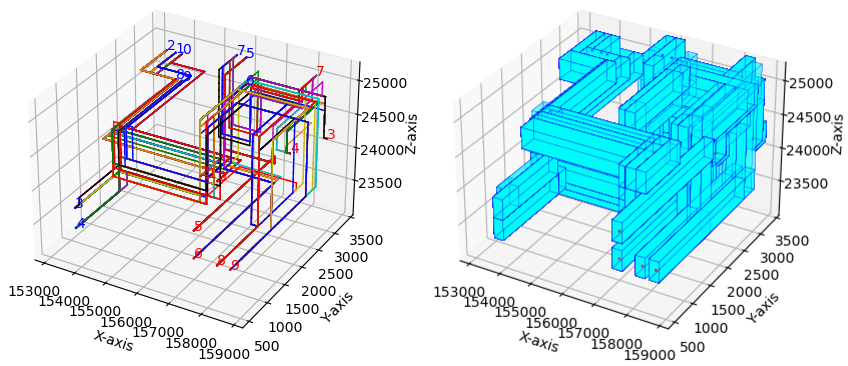}
\caption{ Illustration of the initial region used in the dimensionality reduction algorithm\label{fig:initsol}}
\end{figure}

We propose an iterative procedure that starts by considering some tentative areas as candidate zones for searching for the solution (called in the pseudocode \texttt{Init-Sol}) and then solving the multicommodity flow problem  with additional \eqref{mcf:8} constraints for each single commodity to determine the best solution for each service independently of the remainder services. Then, for each commodity we consider a parallelepiped    of a given initial dimension  around the obtained path and  instead of solving \eqref{prps} in the whole graph, we solve the problem in the union of these parallelepipeds together with the initial candidate zones (see Figure \ref{fig:initsol} (right)).  We do this by fixing  {to zero} the variables indicating arcs outside the above mentioned regions. In case the problem is feasible, it provides a solution of \eqref{prps}. Otherwise, the dimension  of the parallelepipeds is increased and the procedure is repeated until feasibility is obtained. Although in the worst case situation the algorithm requires solving the original instance of the problem (the whole graph), in practice, it allows to solve the problem with a { considerably} smaller number of variables.  Algorithm \ref{alg:reduce} in \ref{appendix:dm} provides a pseudocode for this procedure.

\subsection{Decomposition-based algorithms\label{ss:db}}

Note that when removing the capacity constraints, \eqref{mcf:1}, from MCMNFP, the problem can be decomposed in $|K|$ single-commodity network flow problems with unitary demands which are equivalent to solve $|K|$ shortest path problems on the undirected graph $\widetilde{G}$, and then, solvable in $\mathcal{O}(|V|^2+|E|)$ time by the classical Dijkstra algorithm~\cite{dijkstra1959}.  Based on this idea, we propose an  algorithm that exploits the use of this decomposition for solving fast and accurately the PRPS. Algorithm \ref{alg:mh} in \ref{appendix:db} shows the pseudocode of this  algorithm.

First, the algorithm considers a sorted list of all the commodities (pipelines to route), $K_0$, and a maximum number of iterations to perform, \textit{maxit}. The iterations are grouped in three different classes: \textit{parallel} iterations (\texttt{I\_par}), \textit{sequential} iterations (\texttt{I\_seq}) and \textit{cluster} iterations (\texttt{I\_cluster}), such that the three sets define a partition of $\{1, \ldots, maxit\}$.

The main scheme for all the iterations is similar.
At each iteration $it \in \{1, \ldots, maxit\}$, a shortest path is built for each commodity included in the sorted list $K_{it}$ trying to pass the elbow test, i.e., the minimum required distance between consecutive elbows in the paths (constraints \eqref{mcf:7}). This path is constructed using the procedure \texttt{SPP\_Pipe\_Elbow\_Test} (Algorithm  \ref{alg:codos} in \ref{appendix:db}). In this procedure, for each commodity $k\in K_{it}$, a shortest path from the source node $s^k$ to the destination node $t^k$  is initially built on the undirected graph $\widetilde{G}$ using  Dijkstra algorithm ($\SPP(s^k \rightarrow t^k)$). In case the elbow test is passed, \texttt{SPP\_Pipe\_Elbow\_Test} returns that path. Otherwise, the part of the path that passes the test is kept, the cost of the first  elbow that does not pass the test is increased and  Dijkstra algorithm is run again  (with the new costs) from the last elbow passing the test to the destination node (updating the source node for $\SPP$). The process is repeated every time we find an elbow that does not pass the test. If we identify the process is  cycling we initialize it  maintaining the increased costs and taking the source node of the commodity as destination node and viceversa. This last step is repeated until the elbow test is passed or until a maximum number of iterations (of the elbow test) is reached and then, the elbow test has not been passed.

The second step consists of avoiding overlapping of paths in the solution. This phase is different for each of the types of iterations that we consider.
 In the sequential iterations paths for commodities are solved one by one in the sorted order and therefore, the sorting of the commodities affects the solution. In these iterations, once the path for commodity $k$ is solved,  this path is kept, commodity $k$ is removed from the list of sorted commodities and overlapping is avoided by increasing for the remaining services $k'\neq k$ in the list, the cost of the edges in conflict, that is, the edges  at a distance lower or equal than  $R^k+\Delta^k$ from the edges used in the path of service $k$. 
 They can be tested  using Algorithm \ref{alg:recub} described in  \ref{appendix:db}.

 In the parallel iterations the paths for commodities are  solved independently and the overlapping phase is  addressed once all the paths are constructed by increasing the cost of the edges in conflict.
In the cluster iterations, paths for all the commodities  are solved  initially in parallel.  The commodities that do not overlap with any other are removed from the set of commodities $K_{it}$ and their paths are kept. The remainder commodities $k\in K_{it}$ are organized in clusters of services with edges in conflict and they are sorted by  {some prefixed measure, as for instance,} increasing length of their paths. The first sorted commodity of each cluster is removed for the set of commodities $K_{it}$ and their path is kept. The cost of the edges in conflict with the  paths that are kept is increased for the commodities remaining  in $K_{it}$ and these commodities  are solved in parallel. The process is  repeated until $K_{it}=\emptyset$ or until the maximum number of iterations is reached.

Note that, for all types of iterations, each iteration is performed with a new set of costs  for the edges, trying to avoid the use of conflicting edges in the solutions for the next iteration.

\section{Computational experience}\label{sec:experiments}
Next, we report on the results of some computational experiments that we have run, in order to compare empirically the proposed exact and heuristic approaches.

\subsection{Random instances}
 The set of random instances are generated as follows.
We choose $\mathcal{P}_0$ as a cube of edge length equal to 128 units.
 {An orthogonal 3D grid of $d\times d\times d$ physical nodes, with $d\in\{17,33\}$ is built on this cube, being the distance between adjacent  nodes fixed to $\frac{128}{d-1}$ units. This means a spacing of 8 and 4 units for $d=17$ and $d=33$, respectively. We will call $d$ the \textit{density} of the grid.

Five different groups of 15 obstacles each are generated with
cube shape and edge length $l=10$ units.
Each coordinate of the center of an obstacle is placed in $[\varepsilon+\frac{l}{2},128-\varepsilon-\frac{l}{2}]$ where
$\varepsilon$ stands for the width of an obstacles-free layer along the facets of {$\mathcal{P}_0$} that ensures some feasible origin-destination paths in case too many obstacles are placed. This grid defines a baseline for the undirected graph $\widetilde{G}$  described in Section \ref{sec:PRPND} and  the directed version $G$ described in Section \ref{sec:model}.}

Five different groups of $12$ services are generated with
equal radius of 4 units and
establishing a safety minimum allowed distance of 1 unit with respect to other services.
For each service, a source and a destination have been randomly generated belonging to planes $y=0$ (for the source) and $y=128$ (for the destination)   assuring that both points are nodes of the grid for $d=17$ and therefore, for $d=33$.

Several details on real-life cost functions for our problem are given in Section \ref{sec:case}.
However, for this preliminary study on random instances, we choose the cost function coefficients for an edge $e$ as:
$$
c_{e}^k = \alpha_{1}^k   (d_{e}+  10 El_e+ 2 Ch_e),
$$
where:
\begin{itemize}
  \item $\alpha^{k}_{1}$ takes random integer values in $[1,9]$.
     \item $d_{e}$ is the  physical distance between the two end-nodes $i$ and $j$ of the edge $e=\{i,j\}$ (length of the edge).
   \item $El_e$ takes value 1 if edge $e$ represents an elbow and $0$ otherwise.
    \item $Ch_e$ takes value 1 if edge $e$ is a vertical edge (in the $Z$-axis) and $0$ otherwise.
    \end{itemize}
Observe that the cost function considers the distance of the paths, the number of elbows introduced as well as  changes in height. This cost system is a simplified version of the general cost function  described in Section \ref{sec:case} for our case study.

We denote by $(d,s,o,g)$ the instance of density $d\in\{17,33\}$, services $ \{1,2,...,s\}$ with $s\in\{5,8,12\}$, obstacles $\{1,2,...,o\}$ with $o\in\{5,10,15\}$ and creating five instances, $g\in\{1,\ldots,5\}$, for each combination of services and obstacles. Therefore 90 different benchmark instances are generated.

All instances were solved with the Gurobi 7.7 optimizer, under a Windows 10 environment in an Intel(R) Core(TM)i7 CPU 2.93 GHz processor and 16 GB RAM. Default values were initially used for all parameters of Gurobi solver and a CPU time limit of 7200 seconds was set. We have also tested different combinations of parameters for the solver cut strategy and dimensionality reduction heuristic but, unless it is specified, the best results were obtained with the parameters of the solver set to the default values.
An initial solution was given to the problem by solving the multicommodity flow problem with additional \eqref{mcf:8} constraints for each  service independently of the remainder services.

For the decomposition based heuristic, we fix the maximum number of iterations to $20$, the first $10\%$ of them of type \texttt{I\_par}, the next $80\%$ of type \texttt{I\_cluster} and the last $10\%$ of type \texttt{I\_seq}. In case a feasible solution is not found with the maximum number of iterations (which only happens for some of the $d=33$ instances), we generate the solution for the same instance by reducing the density to $d=17$, which is indeed, a feasible solution of the instance.

In Table \ref{ta:random} we report the average results of our computational experiments.  The first three columns indicate the parameters identifying the instances, $d$ (density), $s$ (number of services) and $o$ (number of obstacles). For each of the combinations of these parameters we report the average results of the five generated instances. Column Vars indicates the  number of variables of the (MCMNFP) problem and column Cons the number of constraints. Column Solved reports the number of instances out of five solved to optimality by the Exact method (Ex). In the block of columns denoted by GAP we report the  gap obtained with the three proposed procedures, the Exact (Ex) model, the dimensionality reduction algorithm (H1), based on trimming down instances, and the one based on decomposing the problem (H2). The GAP reported for the exact approach indicates the MIP gap obtained at the end of the time limit in case the problem has not been optimally solved, while the gap for the heuristic procedures gives the percent deviation of the heuristic solution with respect to the best solution obtained with the exact approach. Finally, in the block of columns Time we report the CPU time required by each one of the approaches.

A first analysis of the results shows that the exact algorithm solves all instances with $d=17$ and any number of services and obstacles, whereas for  $d=33$ we could solve to optimality almost all the instances for number of services $s=5$ and $8$ but not for $s=12$ services. Actually, for $s=12$ services none of the instances could be solved to optimality within the time limit although the final MIP gap is always very small and less than or equal to $1.48\%$.

On the other hand, the results reported for the heuristic algorithms are rather good. For all the instances our two heuristic approaches (H1) and (H2) always find feasible solutions and the gaps with respect to the best solution found by the exact method (Ex) are less than or equal to $5.73\%$. Actually, heuristic (H2) reports even better gaps being always less than or equal to $1.59\%$.

Concerning running times, as expected, the methods based on solving mathematical programming models, i.e., the exact method (Ex) and the matheuristic (H1), require more time to get to their solutions. Clearly, (H1) is less time consuming than (Ex) since it solves trimmed instances which results in programs with much less variables and constraints so that the computing time is also smaller. On the other hand, the decomposition heuristic  (H2) which is based on iteratively solving shortest path problems with modified weights is much lighter and the running times are considerably smaller.


\begin{table}[h!]
\centering
\adjustbox{scale=0.93}{\begin{tabular}{|c|c|c|ccc||ccc||rrr|}
	\hline
	                    &                     &    &                  &                 &             &           \multicolumn{3}{c||}{GAP}           &                        \multicolumn{3}{c|}{Time}                         \\ \cline{7-12}
	         d          &          s          & o  &       Vars       &      Cons       &   Solved    &       Ex       &      H1       &      H2       & \multicolumn{1}{c}{Ex} & \multicolumn{1}{c}{H1} & \multicolumn{1}{c|}{H2} \\ \hline
	\multirow{9}{*}{17} & \multirow{3}{*}{5}  & 5  &      228,232      &      73,572      &      5      &     0.00      &     0.10      &     0.00      &                 11.15 &                   6.26 &                    1.28 \\
	               &                     & 10 &      227,421      &      73,386      &      5      &     0.00      &     0.32      &     0.00      &                 11.51 &                   7.80 &                    1.48 \\
	               &                     & 15 &      226,736      &      73,230      &      5      &     0.00      &     2.01      &     0.00      &                 12.21 &                   9.73 &                    1.40 \\ \cline{2-12}
	               & \multirow{3}{*}{8}  & 5  &      399,406      &     117,761      &      5      &     0.00      &     0.31      &     0.64      &                400.17 &                 396.38 &                    2.24 \\
	               &                     & 10 &      397,986      &     117,463      &      5      &     0.00      &     1.96      &     0.22      &                339.85 &                 240.53 &                    4.37 \\
	               &                     & 15 &      396,788      &     117,214      &      5      &     0.00      &     5.73      &     0.34      &                347.41 &                 482.09 &                    4.59 \\ \cline{2-12}
	               & \multirow{3}{*}{12} & 5  &      627,638      &     176,735      &      5      &     0.00      &     0.61      &     1.59      &               4216.11 &                1824.99 &                    3.78 \\
	               &                     & 10 &      625,407      &     176,289      &      5      &     0.00      &     2.64      &     1.17      &               3367.94 &                2159.86 &                    7.20 \\
	               &                     & 15 &      623,524      &     175,914      &      5      &     0.00      &     4.93      &     1.26      &               2596.96 &                2215.58 &                    7.55 \\ \hline
	\multicolumn{3}{|l|}{\small\textbf{17}}              & \textbf{417,015}  & \textbf{122,396} & \textbf{45} & \textbf{0.00} & \textbf{2.07} & \textbf{0.58} &      \textbf{1225.45} &       \textbf{1149.25} &           \textbf{3.77} \\ \hline
	\multirow{9}{*}{33} & \multirow{3}{*}{5}  & 5  &     1,694,166      &     537,852      &      5      &     0.00      &     0.10      &     0.00      &                163.20 &                   3.31 &                   98.08 \\
	               &                     & 10 &     1,689,282      &     536,568      &      5      &     0.00      &     0.16      &     0.00      &                152.31 &                   5.43 &                   90.72 \\
	               &                     & 15 &     1,684,880      &     535,407      &      5      &     0.00      &     0.69      &     0.00      &                176.34 &                   9.82 &                   98.66 \\ \cline{2-12}
	               & \multirow{3}{*}{8}  & 5  &     2,964,791      &     860,609      &      3      &     0.06      &     0.24      &     0.64      &               5100.74 &                  62.58 &                  183.93 \\
	               &                     & 10 &     2,956,243      &     858,554      &      4      &     0.16      &     1.26      &     0.15      &               4393.18 &                 247.27 &                  181.27 \\
	               &                     & 15 &     2,948,540      &     856,697      &      4      &     0.16      &     3.07      &     0.32      &               3484.85 &                1546.39 &                  184.73 \\ \cline{2-12}
	               & \multirow{3}{*}{12} & 5  &     4,658,958      &     1,291,007     &      0      &     1.26      &     0.45      &     1.59      &               7200.05 &                 236.02 &                  315.66 \\
	               &                     & 10 &     4,645,524      &     1,287,926     &      0      &     1.48      &     1.32      &     1.17      &               7200.05 &                1214.57 &                  311.04 \\
	               &                     & 15 &     4,633,420      &     1,285,139     &      0      &     1.38      &     3.68      &     1.26      &               7200.05 &                4278.22 &                  308.63 \\ \hline
	\multicolumn{3}{|l|}{\small\textbf{33}}              & \textbf{3,097,312} & \textbf{894,418} & \textbf{26} & \textbf{0.50} & \textbf{1.22} & \textbf{0.57} &      \textbf{3896.75} &        \textbf{844.85} &         \textbf{196.97} \\ \hline
	\multicolumn{3}{|c}{\small\textbf{Total}}            & \textbf{1,757,163} & \textbf{508,407} & \textbf{71} & \textbf{0.25} & \textbf{1.64} & \textbf{0.57} &      \textbf{2576.11} &        \textbf{997.05} &         \textbf{100.37} \\ \hline
\end{tabular}}
\caption{Computational results for random instances\label{ta:random}}
\end{table}

\subsection{Case study\label{sec:case}}

This section is devoted to present the application of the proposed methodology to one of the realistic instances provided by our partner Ghenova in order to show the actual applicability of our solution methods to solve scenarios that appear in naval design.

\begin{figure}[h]
\begin{center}
\includegraphics[scale=0.8]{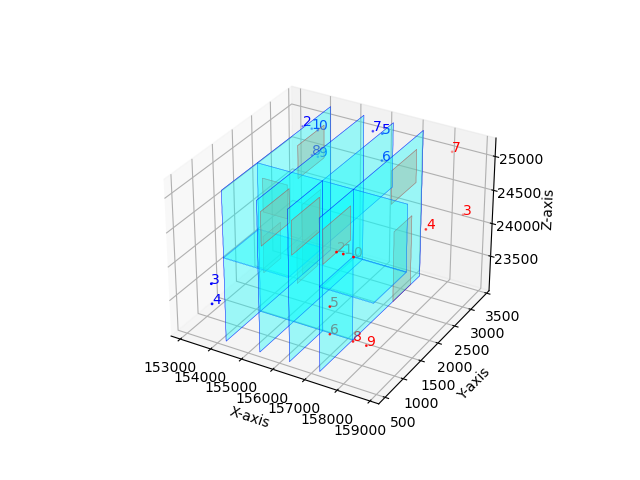}
\caption{Scenario of the case study \label{fig:esc3}}
\end{center}
\end{figure}

In what follows we describe the elements involved in the instance analyzed in this section, which is drawn in Figure \ref{fig:esc3}: 

\begin{itemize}
\item The space to design the pipelines systems consists of the three-dimensional parallelepiped with  widths $5732, 2836$ and $2013$ units  (in the $X$, $Y$ and $Z$ axis, respectively) which represent a ship cabin.
\item A grid was generated for the cabin by subdividing each axis in segments of width $100$, producing an initial grid of $58 \times 29 \times 21$ nodes. The graph $\tilde{G}=(V,E)$ is generated by introducing virtual nodes and edges, with $|V|=73407$ and $|E|=282202$.
\item Ten services are to be routed in such a cabin, each of then identified with a number from $0$ to $9$. The sources  and destinations of the services are located as needed by the designer through the whole cabin (red and blue numbers/services in the figure, respectively).
\item Five obstacles (walls) obstructing the routes are given (light blue shapes in the figure). They consist of metal slices of 10 units width and can be traversed through  $8$ holes/windows (light orange squares in the figure).
\item The minimum distance between consecutive elbows is assumed to be $50$ units, as required by the designer.
\item The radii of the cylindrical pipes are $50$, $75$ and $100$ units and the minimal security distance among two pipelines is $50$ units.
\end{itemize}

As already mentioned, each edge in the generated grid for the instance, has an associated cost that allows to evaluate the feasible routes to be traced in the network. The cost system incorporates the preferences of the designer when routing the pipelines.  We consider an additive cost structure for the objective function to model the use of edges and elbows in the route of each commodity with the following shape:
$$
c_{e}^k = (\alpha_{1}^k + \alpha_5^k Pr_e) d_{e}+ \alpha_2^k  El_e+ \alpha_3^k  H_e+ \alpha_4^k Ch_e+ \alpha_6^k  Pc_e+ \alpha_7^k  Cl_e,
$$
for each edge $e$ and each service $k$. In this function there are some $\alpha$-parameters affected by the characteristics of the edge. The criteria and parameters that define the cost system are detailed in Table \ref{table:costs}. This cost structure is flexible enough to reflect most the preferences of naval designers, and was determined in view of the criteria exposed by our partner company for the selection of implementable solutions.

\begin{table}
\begin{center}
\begin{tabular}{rp{9cm}}\hline
Criteria & Description\\\hline
$d_{e}$ &  physical distance between end-nodes $i$ and $j$ of edge $e=\{i,j\}$ (length of the edge).\\
$El_e$ & 1 if edge $e$ represents an elbow and $0$ otherwise.\\
$H_e=MH-h_e$ &  being $MH$ the maximum height of the ship cabin and $h_e$ the height of edge $e$ in case it is in a plane parallel to the $XY$-plane and $0$ otherwise. \\

$Ch_e$ & 1 if edge $e$ is a vertical edge (in the $Z$-axis) and $0$ otherwise.\\
$Pr_e$ & 1 if edge $e$ belongs to a preference zone and $0$ otherwise.\\
$Pc_e$ & 1 if edge $e$ crosses a penetrable zone and $0$ otherwise.\\
$Cl_e$ & 1 if edge $e$ represents an elbow and it is close to a source or a destination point and $0$ otherwise.\\\hline
Parameters& Description\\\hline
$\alpha_1^k$& Cost per unit length.\\
$\alpha_2^k$& Cost of an elbow.\\
$\alpha_3^k$& Cost of moving away from the ceiling.\\
$\alpha_4^k$& Cost of changing in z-coordinates (height). Therefore, moving to a different height is penalized.\\
$\alpha_5^k$& Bonus per routing the pipeline in a preference zone ($\alpha_5^k < 0$;  $\alpha_1^k + \alpha_5^k > 0$).\\
$\alpha_6^k$& Cost of crossing a penetrable zone.\\
$\alpha_7^k$& Cost of locating an elbow close to the source or destination point of a pipeline.\\\hline
\end{tabular}
\caption{Criteria and parameters involved in the cost function of the case study\label{table:costs}}
\end{center}
\end{table}

\begin{example}\label{ex:cost}
In order to illustrate the above cost function, we use the part of the graph drawn in Figure \ref{figmalla}.
\begin{figure}
\begin{center}
\begin{tikzpicture}[scale=0.8]
\draw[-,black,dashed,line width=0.02cm] (0,3) -- (6,3) -- (10,7) -- (3.85,7) -- (-0.15,3);
\draw[-,black,dashed,line width=0.02cm] (8,8) -- (10,10) -- (4,10) -- (0,6) ;

\draw[-,blue,line width=0.03cm] (5,2) -- (5,11);
\draw[-,red,line width=0.03cm] (2,5) -- (8,5);
\draw[-,black,line width=0.03cm] (4,5) -- (8,5);
\draw[-,red,line width=0.03cm] (2,8) -- (8,8);
\draw[-,red,line width=0.03cm] (3,3) -- (7,7);
\draw[-,red,line width=0.03cm] (3,6) -- (7,10);
\tikzstyle{linea1}=[-,green,line width=0.03cm]
\tikzstyle{nodo1}=[circle, draw=blue!80, fill=white,inner sep=0.5mm]
\tikzstyle{nodo2}=[circle, draw=blue!80, fill=white,inner sep=0.5mm]

\node[right]    () at (8,5) { $\Pi_2$};
\node[right]    () at (8,8) { $\Pi_1$};

\node[nodo2]    () at (8,5) {}; \node[nodo2]    () at (2,5) {};
\node[nodo2]    () at (5,2) {}; 
\node[nodo2]    () at (3,3) {}; \node[nodo2]    () at (7,7) {};

\node[nodo2]    () at (8,8) {}; \node[nodo2]    () at (2,8) {};
\node[nodo2]    () at (5,11) {}; 
\node[nodo2]    () at (3,6) {}; \node[nodo2]    () at (7,10) {};

\node[nodo1]    (x2) at (4,5) {};
\node[nodo1]    (y1) at (5,6) {};
\node[nodo1]    (z1) at (4.75,4.75) {};
\node[nodo1]    (x22) at (4,8) {};
\node[nodo1]    (y12) at (5,9) {};
\node[nodo1]    (z12) at (4.75,7.75) {};
\draw[-,linea1] (x2) -- (y1);
\draw[-,linea1] (x2) -- (z1);
\draw[-,linea1] (y1) -- (z1);
\draw[-,linea1] (x22) -- (y12);
\draw[-,linea1] (x22) -- (z12);
\draw[-,linea1] (y12) -- (z12);

\node[above right] () at (x2) {$v_2$};
\node[above right] () at (x22) {$v_1$};
\end{tikzpicture}
\end{center}
\caption{Graph used to illustrate the cost function in Example \ref{ex:cost} \label{figmalla}}
\end{figure}
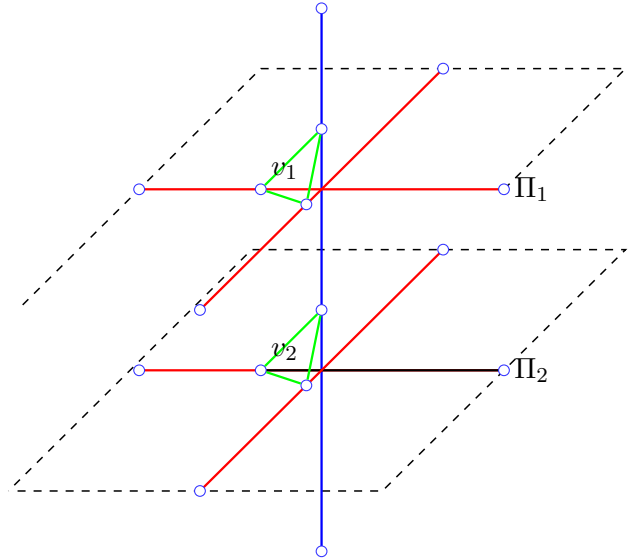

In this graph,  $\Pi_1$ is the plane $\{(x,y,z): z=1\}$ while $\Pi_2 =\{(x,y,z): z=0\}$ and then, plane $\Pi_1$ is closer to the ceiling than plane $\Pi_2$. We denote by $v_1$ and $v_2$ the nodes which are virtualized in the figure (for each of the hyperplanes). We assume that there is a preference zone using $\Pi_1$ but there are no preference zones using $\Pi_2$. Furthermore, $v_2$ is assumed to be close to  the source or  destination point of  some service and also that the black edge in $\Pi_2$ crosses a penetrable zone for the service. In this situation, we have:

\begin{itemize}
\item The criteria $d$ and $El$  coincide for all the edges that belong to the planes $\Pi_1$ and $\Pi_2$, being the part of their costs affected by parameters $\alpha_1$ and $\alpha_2$, the same.
\item The height of the edges in $\Pi_1$ is larger than height of the edges in $\Pi_2$, being the term $\alpha_3 H$ smaller for the edges in $\Pi_1$.
\item The criterion $Pr$ affects the edges in $\Pi_1$ but not those in $\Pi_2$.
\item The criteria $Pc$ and $Cl$ affect  some of the edges in $\Pi_2$ but none of  those in $\Pi_1$.
\item The edges linking nodes in $\Pi_1$ with nodes in $\Pi_2$ (parallel to the $Z$-axis) are only affected by criteria $d$ and $Ch$.
\end{itemize}
In Table \ref{table:cost} we summarize the cost function for edges in the different planes depicted in Figure \ref{figmalla}.
\begin{table}
\begin{center}
\begin{tabular}{|c|cc|}\cline{2-3}
\multicolumn{1}{c|}{}& Virtual Arcs (Elbows) & Non Virtual Arcs (Straight Pipes)\\\hline
$e\in \Pi_1$ & $\alpha_2^k$ & $(\alpha_1^k+\alpha_5^k) d_{e} + \alpha_3^k H_e$\\\hline
 \multirow{2}{*}{$e\in \Pi_2$} & \multirow{2}{*}{$\alpha_2^k + \alpha_7^k$} & $\alpha_1^k d_{e} + \alpha_3^k H_{e} $ \\
& & $\alpha_1^k d_{e} + \alpha_3^k H_{e}+ \alpha_6^k$ (penetrable edges)\\\hline
e parallel to $Z$-axis & & $\alpha_1^k d_{e} + \alpha_4^k$\\\hline
\end{tabular}
\caption{Illustration of the cost function\label{table:cost}}
\end{center}
\end{table}

\end{example}

In the actual instance considered in the case study, each edge in the graph $ \widetilde{G}$ has the same cost for all services, i.e., $c_{e}^k = c_{e}^{k^\prime}$, for all $k, k^\prime \in K$, as required by our partner company. The parameters provided for this instance were:
\begin{center}
\begin{tabular}{|ccccccc|}\hline
$\alpha_{1}$ & $\alpha_{2}$  & $\alpha_{3}$ & $\alpha_{4}$  & $\alpha_{5}$  & $\alpha_{6}$  & $\alpha_{7}$  \\\hline
1                            & 2800                        & 700                          & 200                          & -0.7                         & 4000                         & 3000\\\hline
\end{tabular}
\end{center}
and the preference zones were defined as the following parallelepipeds:
$
[154241,158844] \times [535,1419] \times [24249, 25221] \text{ and  } [152013,154241] \times [535,3371] \times [24249, 25221].
$

%

We run our model for this instance with the exact and the different heuristic approaches. A time limit of 12 hours was fixed for all the procedures. The exact branch-and-cut approach was not able to find a feasible solution of the problem within the time limit. Note that the number of variables for this problem is $2819480$ and the (initial) number of constraints is $31014280$ (constraints \eqref{mcf:6} and \eqref{mcf:8} that are added as lazy constraints in the branch-and-cut approach are not counted here). The number of explored nodes after 12 hours was $1804$, which gives an idea of the high computational load required to solve the problem at each node of the branch-and-bound tree.

For the dimensionality reduction matheuristic (H1), we first 
solve, independently, the multicommodity flow problem with additional \eqref{mcf:8} constraints for each single service, and restrict the search region of our problem to the region induced by these paths enlarged by a parallelepiped of dimension $\delta_k=\max\{R^{k'}+\Delta^{k'}: k' \in K\}$ around the paths. The value of $\delta_k$ is sequentially enlarged by $0.1$ units until a feasible solution is obtained. This procedure obtains a solution in $132.37$ seconds, and the last problem explored $18041$ nodes in the branch-and-bound tree. The number of variables of this problem was $473292$ ($16.78\%$ of the variables of the exact approach) and the initial number of constraints was $807657$ ($2.6\%$ of the number of constraints required in the exact approach). The number of lazy constraints added in the branch-and-cut approach was $40800$ ($12$ of them of type \eqref{mcf:8}). The objective value of the solution was $448690.78$.

For the decomposition based heuristic (H2) we fix to 10 the maximum number of iterations with the first $10\%$ of them of type \texttt{I\_par}, the next $20\%$ of type \texttt{I\_cluster} and the last $70\%$ of type \texttt{I\_seq}. This configuration was adequately tuned using a simplified instance for the problem as a training sample. This heuristic computed a solution in $94.54$ seconds after $6$ iterations and the obtained objective value was $433387.86$. Observe that the deviation between the solutions obtained with the two heuristics was $3.41\%$.

\begin{figure}
\begin{center}
\includegraphics[scale = 0.7]{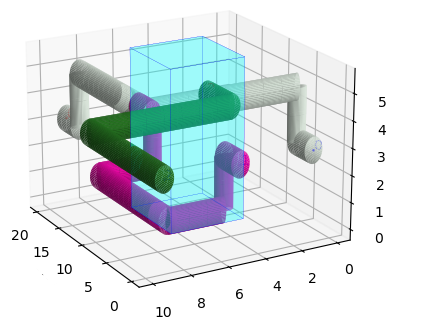}~\hspace*{-0.03cm}
\includegraphics[scale = 0.7]{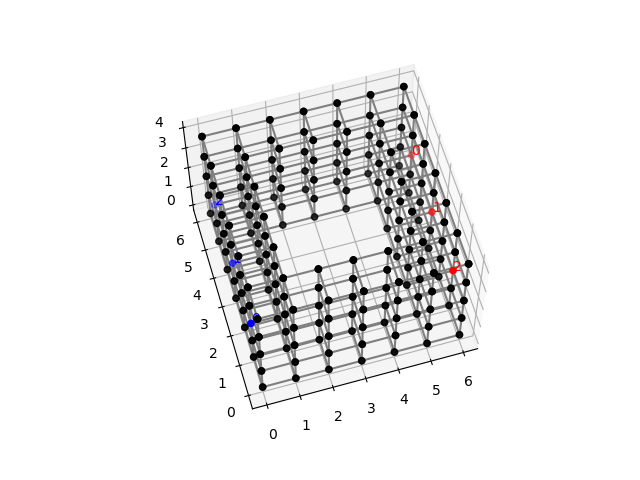}
\caption{Graphical display of the solution for the case study with 10 services.  In the left it is shown the solution without showing the obstacles. The figure on the right shows the same solution integrated within the set of obstacles in the scenario.\label{fig:1sol-esc3}}
\end{center}
\end{figure}

\begin{figure}
\begin{center}
\includegraphics[scale = 0.75]{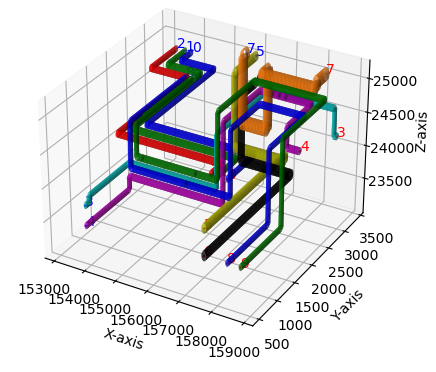}~
\includegraphics[scale = 0.75]{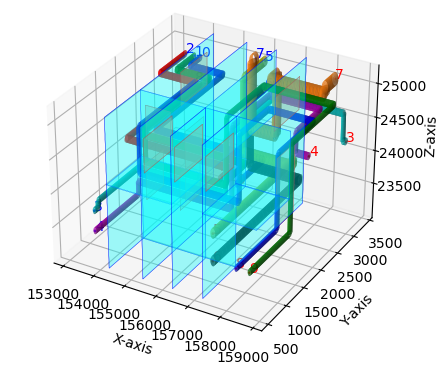}
\end{center}
\caption{Graphical display of the solution for the case study with 10 services actual pipelines.  In the left it is shown the set of pipelines without showing the obstacles. The figure on the right shows the same solution integrated within the set of obstacles in the scenario.\label{fig:2sol-esc3}}
\end{figure}

We draw in Figure \ref{fig:1sol-esc3} the obtained solution (with heuristic H1) as paths in the graph. In the picture we represent the obtained paths without (left) and with  obstacles and holes (right). The obtained solution as actual pipelines is drawn in Figure \ref{fig:2sol-esc3}. One can observe from the picture that real-world instances require routing intricate pipelines avoiding obstacles and respecting the space between the pipes what really makes difficult to locate them in the cabin. The heuristic approaches produce, in view of the results obtained with the synthetic instances, good quality solutions in reasonable CPU times, providing the naval designer with a powerful decision-aid  tool for this task.

\section{Conclusions}\label{sec:conclusions}

This paper develops a general methodology for the optimal automatic network design of different commodities incorporating real-world constructability conditions. This methodology was originated by an actual  problem consisting of the design of  spatial pipelines with ships cabins motivated by a recent collaboration with the company  Ghenova, a leading naval engineering company. Our proposal is to adapt an \textit{`ad hoc'} minimum cost multicommodity network flow model to the problem that includes all the technical requirements of feasible pipelines routing. The large number of variables and constraints of real-sized instances makes it inappropriate  (impossible in real-world instances) to load into the solver the complete model so that we implemented a branch-and-cut algorithm, where, initially only the standard multicommodity flow model is considered and additional technical constraints are separated within the branch-and-bound tree as lazy constraints. On top of that, we also develop two heuristic algorithms that provide rather good feasible solutions. The first one is a matheuristic based on solving trimmed instances and the second one is a decomposition heuristic that iteratively solves shortest path problems with modified edges weights. Our computational results on randomly generated instances are rather promising showing that we can solve to optimality medium-sized instances. We also report a case study based on a realistic naval instance of a ship cabin provided by our partner company to our research group to test the methodology.

Future work on this topic includes the consideration of more general underline networks including non-orthogonal designs which are conceptually easy to handle once the admissible solutions graphs are generated whereas the generation of these graphs is still a challenging question that will be the focus of a follow up paper.

We will also explore the use of Machine Learning tools to determine, in advance, the parameters of the cost function, that are assumed in this paper to be given by the designer. The study of the subregions of $\alpha$-parameters inducing the same routes will allow us to provide the designers different options of reasonable parameters to explore in the decision making process.


\section*{Acknowledgements}
The authors of this research acknowledge financial support by the Spanish Ministerio de Ciencia y Tecnologia, Agencia Estatal de Investigacion and Fondos Europeos de Desarrollo Regional (FEDER) via project PID2020-114594GB-C21. The authors also acknowledge partial support from projects FEDER-US-1256951, Junta de Andalucía P18-FR-422, CEI-3-FQM331, NetmeetData: Ayudas Fundación BBVA a equipos de investigación cient\'ifica
2019, and Contrataci\'on de Personal Investigador Doctor (Convocatoria 2019) 43 Contratos Capital Humano L\'inea 2. Paidi 2020, supported by the European Social Fund and Junta de Andaluc\'ia.


\newpage

\appendix

\section{Pseudocode of the algorithms described in Section\ref{sec:matheuristic}}\label{appendix}
We gather in this section the pseudocode of the different algorithms described in Section \ref{sec:matheuristic} devoted to our  matheuristic algorithms.
\subsection{Pseudocode for the dimensionality reduction  algorithm in Section \ref{ss:dm}}\label{appendix:dm}
Next, we include the pseudocode of the dimensionality reduction algorithm.
 \begin{algorithm}[h]
\SetKwInOut{Input}{input}\SetKwInOut{Output}{output}
\LinesNumbered
\Input{
Initial set of candidate zones provided by user $\{{\rm \texttt{Init-sol}}_k\}_{k\in K}. $\\
$\{\delta_k=R^k+\Delta^k\}_{k\in K}$ initial dimension of the parallelepiped around the best  path for commodity $k$. \\
it$=$1.
$stop=0$.
}

\For{$k \in K$}{
Solve the Multicommodity Flow Problem with additional \eqref{mcf:8} constraints for the single commodity $k$ in the whole graph $G$: $\{{\rm path}_k\}_{k\in K}$\\
}

\While{$stop=0$}{

Solve \eqref{prps}(it), where all variables $x_{ij}^k$ for a given service $k$, that are outside ${\rm Init-sol}_k$ union with the parallelepiped of dimension  $\delta_k$ around ${\rm path}_k$, are fixed to zero.

\eIf{\eqref{prps}{\rm (it)} is feasible}{
$stop=1$
}{
Increase $\delta_k$ for all $k$ in $K$
}
it$=$it+1

}

\Output{Solution of \eqref{prps}(it)}
\caption{Dimensionality reduction algorithm.\label{alg:reduce}}
\end{algorithm}
\newpage

\subsection{Pseudocode for the decomposition-based algorithm described in  Section \ref{ss:db}}\label{appendix:db}
\noindent In the following we describe the different modules of the decomposition-based  matheuristic.

\begin{algorithm}
\SetKwInOut{Input}{input}\SetKwInOut{Output}{output}
\LinesNumbered
\Input{$stop:=0$ (termination criterion); \texttt{it}=0 (iterations counter); $\{1, \ldots, maxit\}$ = $\texttt{I\_par} \cup \texttt{I\_seq} \cup \texttt{I\_cluster}$ (Index ser for iterations); $K_0$ (sorted set of all services in $K$).}

\While{$stop=0$ and $\texttt{it}\leq$ maxit}{
    \If{$\texttt{it} \in \texttt{I\_seq}$ }{
        $K_{it}=K_0$ and $\bar K=K_0$.
    }
    \For{$k \in K_{it}$}{
        Apply the \texttt{SSP\_Pipe\_Elbow\_Test} to $k$. (Algorithm \ref{alg:codos})\\
        Update \texttt{Covering\_List($k$)}. (Algorithm \ref{alg:recub}) \\
        \If{$\texttt{it} \in \texttt{I\_seq}$}{
            $\bar K=\bar K\setminus \{k\}$.
            For all $k^\prime  \in  \bar K$, slightly increase the cost of edges linking vertices in  \texttt{Covering\_List($k$)}.
        }
    }
    \For {$k \in K_{it}$ and $k^\prime \in K$}  {
        Compute  $\texttt{Cov}(k,k^\prime):=$\texttt{Covering\_List($k$)}$\cap$\texttt{Covering\_List($k^\prime$)}.
    }
    \If { \rm $\texttt{Cov}(k,k^\prime)\neq \emptyset$ for some pair $(k,k^\prime)$,  $k \in K_{it}$, $k^\prime \in K$}{
        $stop=0$.\\
        \If {$\texttt{it} \in \texttt{I\_par}$}{
            For all service  $k \in K_{it}$,
                increase the cost (to the most costly service) of edges linking vertices in $\bigcup_{k^\prime\in K}\texttt{Cov}(k,k^\prime)$.
        }
        \If {\rm $\texttt{it} \in \texttt{I\_cluster}$}{
            $K_{it+1}=\emptyset;$ \\
            Sort the services, $k_1 \succ k_2 \succ \cdots \succ k_\ell$ (by priority);\\
            \If {\rm $\texttt{Cov}(k_i,k_j)\neq\emptyset$  and $k_i \succ k_j$}{
                $K_{it+1}=K_{it+1} \cup\{k_j\}$.\\
            }
            For all service $k\in K_{it+1}$,
                increase the cost (to the most costly service) of edges linking vertices in $\bigcup_{k^\prime\in K\setminus K_{it+1}}\texttt{Cov}(k,k^\prime)$.
        }
        $\texttt{it}=\texttt{it}+1$
    }
    \Else{
        $stop=1$
    }
}

\caption{Decomposition-based Algorithm.\label{alg:mh}}
\end{algorithm}

\begin{algorithm}
\SetKwInOut{Input}{input}\SetKwInOut{Output}{output}
\LinesNumbered
\Input{$N^0:=s^k$ (source node for service $k$); $N^f:=t^k$ (destination node for service $k$); $P_k=\emptyset$ (initial path for service  $k$); $test:=1$ (initial  {in}feasibility test).
}

\While{$(test=1)$}  {
$P_k:=P_k\cup \SPP(N^0\rightarrow N^f)$;\\
\If{\rm \eqref{mcf:7} is verified by $P_k$}{
$test:=0$ and
\texttt{return} $P_k$}
\Else{
Identify previous node that has passed the elbow test: $i_0$;\\
$N^0:=i_{0}$;\\
$P_k:=P_k\backslash\SPP(i_{0}\rightarrow N^f)$;\\
Increase the cost of the elbow not verifying the elbow test;
}
}

\caption{\texttt{SPP\_Pipe\_Elbow\_Test}.\label{alg:codos}}
\end{algorithm}

\begin{algorithm}
\SetKwInOut{Input}{input}\SetKwInOut{Output}{output}
\LinesNumbered
\SetKw{KwGoTo}{go to}

\Input{$P_k$= path connecting source node $s^k$ and destination node $t^k$; $L_k=P_k$ and $\delta_k=R^k+\Delta^k$ (covering of the pipeline).}
\For{$v_0 \in P_k$}{
    \For{$v \in L_k$\label{etiq}}{
        \While{ there exists a vertex  $v^\prime$ adjacent to $v$ such that $d_{v^\prime v_0}<\delta_k$ and $v^\prime \notin L_k$ }{
            $L_k=L_k\cup\{v^\prime\}$;\\
             \textbf{go to line} 2.
        }
    }
}
\Output{$L_k$.}
\caption{\texttt{Covering\_List($k$)}. \label{alg:recub}}
\end{algorithm}


\begin{thebibliography}{}

\bibitem[Ando and Kimura, 2011]{and11}
Ando, Y. and Kimura, H. (2011).
\newblock An automatic piping algorithm including elbows and bends.
\newblock In {\em International Conference on Computer Applications in
  Shipbuilding (ICCAS)}, pages 153--158, Triestre, Italy.

\bibitem[Asmara, 2013]{asm13}
Asmara, A. (2013).
\newblock {\em Pipe routing framework for detailed ship design}.
\newblock PhD thesis, Delft University of Technology, Netherlands.

\bibitem[Asmara and Nienhuis, 2006]{asm06}
Asmara, A. and Nienhuis, U. (2006).
\newblock Automatic piping system in ship.
\newblock In {\em International Conference on Computer and IT Applications
  (COMPIT)}, Delft, Netherlands.

\bibitem[Cuervas et~al., 2017]{Cuervas}
Cuervas, F., Tordera, A., Font{\'a}n, A., Brenes, P., Alejo, C., Tovar, A.,
  Puerto, J., Conde, E., Ortega, F., and Hinojosa, Y. (2017).
\newblock Ariadna: sistema autom{\'a}tico de trazado de tuber{\'\i}as y
  canalizaciones en ingenier{\'\i}a.
\newblock {\em Ingenier{\'\i}a Naval}, 963:75--84.

\bibitem[Dijkstra, 1959]{dijkstra1959}
Dijkstra, E.~W. (1959).
\newblock A note on two problems in connexion with graphs.
\newblock {\em Numerische mathematik}, 1(1):269--271.

\bibitem[Even et~al., 1975]{even1975complexity}
Even, S., Itai, A., and Shamir, A. (1975).
\newblock On the complexity of time table and multi-commodity flow problems.
\newblock In {\em 16th Annual Symposium on Foundations of Computer Science
  (sfcs 1975)}, pages 184--193. IEEE.

\bibitem[Guiradello and Swaney, 2005]{gui05}
Guiradello, R. and Swaney, R.~S. (2005).
\newblock Optimization of process plant layout with pipe routing.
\newblock {\em Computers and Chemical Engineering}, 30:99--114.

\bibitem[Hightower, 1969]{H69}
Hightower, D. (1969).
\newblock A solution to line routing problems on the continuous plane.
\newblock In {\em Proceedings of Sixth Annual Design Automation Conference.
  IEEE}, pages 1--24.

\bibitem[Ikehira et~al., 2005]{IKK05}
Ikehira, S., Kimura, H., and Kajiwara, H. (2005).
\newblock Automatic design for pipe arrangement using multi-objective genetic
  algorithms.
\newblock In {\em International Conference on Computer Applications in
  Shipbuilding (ICCAS)}, pages 97--110, Busan, Korea, 23-26 August.

\bibitem[Ito, 1999]{Ito99}
Ito, T. (1999).
\newblock A genetic algorithm approach to piping route planning.
\newblock {\em Journal of Intelligence Manufacturing}, pages 103--114.

\bibitem[Kim et~al., 2013]{kim13}
Kim, S.~H., Ruy, W.~S., and Jang, B.~S. (2013).
\newblock The development of a practical pipe auto-routing system in a
  shipbuilding cad environment using network optimization.
\newblock {\em Int. J. Naval Archit. Ocean Eng.}, 5:468--477.

\bibitem[Kimura, 2011]{Kim11}
Kimura, H. (2011).
\newblock Automatic designing system for piping and instruments arrangement
  including branches of pipes.
\newblock In {\em International Conference on Computer Applications in
  Shipbuilding (ICCAS2011)}, volume~3, pages 93--99.

\bibitem[Lee, 1961]{lee61}
Lee, C.~Y. (1961).
\newblock An algorithm for path connections and its applications.
\newblock {\em IEEE Transactions on Electronic Computers}, 10(3):346--365.

\bibitem[Min et~al., 2020]{Min2020}
Min, J., Ruy, W., and Park, C. (2020).
\newblock Faster pipe auto-routing using improved jump point search.
\newblock {\em International Journal of Naval Architecture and Ocean
  Engineering}, pages 596--604.

\bibitem[Mitsuta et~al., 1986]{mit86}
Mitsuta, T., Kobayashi, Y., Wada, Y.and~Kiguchi, T., and Yoshinaga, T. (1986).
\newblock A knowledge-based approach to routing problems in industrial plant
  design.
\newblock In {\em Proceedings of the Sixth International Workshop Expert System
  and Their Applications}, pages 237--256, Avignon, France, March.

\bibitem[Newell, 1973]{new73}
Newell, R.~G. (1973).
\newblock {\em Algorithms for the design of chemical plant layouut and pipe
  routing}.
\newblock PhD thesis, Imperial College, London, England.

\bibitem[Park, 2002]{P02}
Park, J. (2002).
\newblock Pipe-routing algorithm development for a ship engine room design.
\newblock In {\em Ph.D. Washington University}.

\bibitem[Park and Storch, 2002]{par02}
Park, J.~H. and Storch, R.~L. (2002).
\newblock Pipe-routing algorithm development: case study of a ship engine room
  design.
\newblock {\em Expert Systems with Applications}, 23:299--309.

\bibitem[Qian et~al., 2008]{QRW08}
Qian, X., Ren, T., and Wang, C. (2008).
\newblock A survey of pipe routing design.
\newblock In {\em Control and Decision Conference, Chinese IEEE}, pages
  3994--3998, Yantai, Shandong, China, 2-4 July.

\bibitem[Rourke, 1975]{rou75}
Rourke, P.~W. (1975).
\newblock {\em Development of a three-dimensional pipe routing algorithm}.
\newblock PhD thesis, Lehigh University.

\bibitem[Salimifard and Bigharaz, 2020]{Sal2020}
Salimifard, K. and Bigharaz, S. (2020).
\newblock The multicommodity network flow problem: state of the art
  classification, applications, and solution methods.
\newblock {\em Operational Research}.

\bibitem[Schmidt-Traub et~al., 1999]{sch99}
Schmidt-Traub, H., Holtk{\"o}tter, T., Lederhose, M., and Leuders, P. (1999).
\newblock An approach to plant layout optimization.
\newblock {\em Chemical Engineering Technology}, 22(2):105--109.

\bibitem[Shing and Hu, 1986]{shi86}
Shing, M.~T. and Hu, T.~C. (1986).
\newblock Computational complexity of layout problems.
\newblock In {\em Layout design and verification}, pages 267--294. Elsevier
  Science Publishers B.V., North Holland.

\bibitem[Singh and Cheng, 2021]{Singh2021}
Singh, J. and Cheng, J. (2021).
\newblock Automating the generation of 3d multiple pipe layout design using bim
  and heuristic search methods.
\newblock In {\em 18th International Conference on Computing in Civil and
  Building Engineering. ICCCBE 2020}, volume~98, pages 54--72.

\bibitem[Xiaoning et~al., 2006]{XYZ06}
Xiaoning, F., Yan, L., and Zhuoshang, J. (2006).
\newblock The ant colony optimization for ship pipe route design in 3d space.
\newblock In {\em World Congress on Intelligent Control and Automation
  (WCICA)}, pages 3103--3108, Dalian, China, 21-23 June.

\bibitem[Xiaoning et~al., 2007]{XYZ07}
Xiaoning, F., Yan, L., and Zhuoshang, J. (2007).
\newblock Ship pipe routing design using the aco with iterative pheromone
  updating.
\newblock {\em Journal of Ship Production}, 23(1):36--45.

\bibitem[Zhu and Latombe, 1991]{zhu91}
Zhu, D. and Latombe, J. (1991).
\newblock Pipe routing-path planning (with many constraints).
\newblock In {\em IEEE International Conference on Robotics and Automation
  (ICCAS)}, pages 1940--1947, Sacramento, California.

\end{thebibliography}
\end{document}